\documentclass{elsart}
\usepackage{amssymb}
\usepackage{amsmath}
\usepackage[mathscr]{eucal}
\usepackage{enumerate}
\usepackage{latexsym}

\newcommand{\R}{\mathbb{R}}

\newcommand{\N}{\mathbb{N}}

\newcommand{\Z}{\mathbb{Z}}
\newcommand{\Q}{\mathbb{Q}}

\newcommand{\M}{\mathbb{M}}
\newcommand{\X}{\mathbb{X}}
\newcommand{\D}{\mathbb{D}}
\newcommand{\skel}{\mathrm{skel}}
\newcommand{\age}{\mathrm{age}}

\newtheorem{definition}{{\bf Definition}}
\newtheorem{theorem}{{\bf Theorem}}
\newtheorem{proposition}{\noindent {\bf Proposition}}
\newtheorem{corollary}{\noindent {\bf Corollary}}

\newtheorem{problem}[definition]{\noindent{\bf Problem}}
\newtheorem{problems}[definition]{\noindent{\bf Problems}}

\newtheorem{lemma}[definition]{\noindent {\bf Lemma}}
\def\endproof{\hfill {\kern 6pt\penalty 500
\raise -0pt\hbox{\vrule \vbox to5pt {\hrule width 5pt
\vfill\hrule}\vrule}}}

\begin{document}
\begin{frontmatter}
\title{Representation of ideals of relational structures}
\author{Christian Delhomm\'e}
\address{E.R.M.I.T., D\'epartement de Math\'ematiques et d'Informatique, Universit\'e de La R\'eunion, 15, avenue Ren\'e Cassin, BP 71551,  97715 Saint-Denis Messag. Cedex~9, France}
\ead{delhomme@univ-reunion.fr}
\author{Maurice Pouzet\thanksref{label2}}
\address{PCS, Universit\'e Claude-Bernard Lyon1,
Domaine de Gerland -b\^at. Recherche [B], 50 avenue Tony-Garnier, F69365  Lyon cedex 07, France, }
\thanks[label2]{Supported by Intas. Research done while the author visited the Math. dept. of U of C in spring 2005 and summer 2006; the support provided is gratefully aknowledged.}
 \ead{ pouzet@univ-lyon1.fr }

\author{G\'abor S\'agi\thanksref{label4}}
\address{Alfr\'ed R\'enyi Institute of Mathematics, 
Hungarian Academy of Sciences,
Budapest H-1364, pf. 127, Hungary}
\thanks[label4]{Supported by Hungarian National Foundation for
Scientific Research grants D042177, T35192.}

\ead{sagi@renyi.hu}
\author{Norbert Sauer\corauthref{cor}\thanksref{label3}}
\address{Department of Mathematics and Statistics, The University of Calgary, Calgary, T2N1N4, Alberta, Canada}
 \thanks[label3]{Supported by NSERC of Canada Grant \# 691325}
\corauth[cor]{Corresponding author.}
\ead{nsauer@math.ucalgary.ca}

\date{\today}

\begin{abstract} The \textit{age} of a relational structure  $\mathfrak A$  of signature $\mu$ is the set
$age(\mathfrak A)$ of its finite induced substructures,
considered up to isomorphism. This is an ideal in the poset
$\Omega_\mu$ consisting of finite structures of  signature $\mu$ and
ordered by embeddability. If the structures are made of infinitely
many relations and if, among those,  infinitely many are at least
binary then there are ideals which do not come from an age. We provide many examples. We particularly look at  metric spaces and offer several problems. We also provide an example
of an ideal $I$ of isomorphism types of at most countable
structures whose signature consists of a single ternary relation
symbol.  This ideal   does not come from the set $\age_{\mathfrak
I }(\mathfrak A)$ of isomorphism types of substructures of
$\mathfrak A$ induced on the members of an ideal $\mathfrak I$ of
sets. This answers a question due to R.~ Cusin and J.F. ~Pabion
(1970).
\end{abstract}

\end{frontmatter}

\section{Introduction and  basic notions}

Let $\mathbb{N}$ be the set of non-negative integers, $\N^*:= \N\setminus \{0\}$ be the set of positive integers and $n\in \N^*$.  A {\em $n$-ary relation} on a set $A$ is a map $R: A^n\to  \{0,1\}$. A {\em
signature}  is a function $\mu:I\to \mathbb{N^*}$ from an {\em index
set } $I$ into $\N$. We write $\mu=(\mu_i;
i\in I)$ as an indexed set.  A {\em relational structure with signature
$\mu$}  is a pair  $\mathfrak{A}:=(A; \mathbf{R}^\mathfrak{A})$
where $\mathbf{R}^{\mathfrak{A}}:=
(\mathsf{R}^{\mathfrak{A}}_i)_{i\in I}$  is a set of  relations on
$A$, each relation $\mathsf{R}^\mathfrak{A}_i $ having arity
$\mu_i$.  If $\mathfrak{A}$ is clear from the context then we will
write $\mathbf{R}$ instead of $\mathbf{R}^\mathfrak{A}$ and
$\mathsf{R}_i$ instead of $\mathsf{R}^\mathfrak{A}_i$. As much as
possible we will denote relational structures by letters of the
form $\mathfrak{A}$, $\mathfrak{B}$, $\mathfrak{C}$, etc. and the
corresponding base sets by $A, B, C$, etc. The {\em cardinality }
of the relational structure $\mathfrak{A}:=(A; \mathbf{R})$ is the
cardinality of $A$, which, as usual, will be denoted by $|A|$.

The signature $\mu=(\mu_i; i\in I)$ is {\em unary, binary,
ternary} and in general {\em $n$-ary} if the range of the function
$\mu$ is $\{1\}$, $\{2\}$, $\{3\}$ or in general $\{n\}$. The
signature $\mu=(\mu_i; i\in I)$ is at most {\em  binary, ternary}
and in general {\em $n$-ary} if the range of the function $\mu$ is
a subset of $\{1,2\}$, $\{1,2,3\}$ and in general $\{1,2, 3,
\dots, n\}$. The signature $\mu$ is \emph{finite} if the index set $I$ is
finite, it is {\it infinite} otherwise. If $S\subseteq \mathbb{N^*}$ then $\mu^{-1}[S]$ is the set
of all indices $i\in I$ for which $\mu_i\in S$. If
$|\mu^{-1}(\mathbb{N^*})|=1$ then $\mu$ is a \emph{singleton signature}.
For example, a relational structure with  a binary singleton
signature is a directed graph which may have loops.

A relational structure is {\em unary, binary, ternary, $n$-ary}
and so on if its signature is unary, binary, ternary, $n$-ary,
respectively.

Let $\mathfrak{A}:=(A; \mathbf{R}^\mathfrak{A})$ and
$\mathfrak{B}:=(B;\mathbf{R}^\mathfrak{B})$ be two relational
structures with common signature $\mu:=(\mu_i;i\in I)$. A map 
 $f:A\to B$ is an {\em isomorphism } of $\mathfrak{A}$ onto $\mathfrak{B}$ if $f$ is bijective  and for all $i\in I$ and
$(x_0,x_1,\dots,x_{\mu_{i}-1})\in A^n$:
\[
R^\mathfrak{A}_i(x_0,x_1,\dots,x_{\mu_{i}-1}) \text{\, \, \,  if
and only if \, \, \,  } R^\mathfrak{B}_i(f(x_0,f(x_1),\dots,
f(x_{\mu_{i}-1})).
\]
Let $A'$ be a subset of $B$, the {\em substructure of
$\mathfrak{B}$   induced on  $A'$},  also called the {\it
restriction of $\mathfrak{B}$ to $A'$} is the relational structure  $\mathfrak{ B}_{\restriction A' }:= (A';  \mathbf{R}^\mathfrak{A'})$, where $\mathbf{R}^\mathfrak{A'}_{i}$ is the restriction of the map $\mathbf{R}^\mathfrak{B}_i$ to $A'^{n_{i}}$.  A map 
 $f:A\to B$ is an {\em embedding  } of $\mathfrak{A}$ into $\mathfrak{B}$ if $f$ is an isomorphism from $\mathfrak A$ onto $\mathfrak B_{\restriction f(A)}$. 
We write $\mathfrak{A}\leq \mathfrak{B}$ to indicate that there exists an embedding of $\mathfrak{A}$ into $\mathfrak{B}$. 

Two relational structures  are {\em isomorphic} or have the same {\em isomorphism type} if there is an isomorphism of one onto the other. We suppose that   isomorphism  types have been defined, we will denote $Is({\mathfrak A})$ denotes the isomorphism type of ${\mathfrak A}$
and we will denote by $\Omega^\ast_\mu$  the class of isomorphism types of
relational structures with signature $\mu$. We will denote by $\Omega_\mu$ the subclass made of isomorphism types of finite relational structures. This class turns out to be a set (of size $\aleph_0$ if $I$ is finite, and of size $2^{\vert I\vert }$ otherwise).  
The relation $\leq$ is a quasi-order on the class of
relational structures with  signature $\mu$. It induces a quasi-ordering on the class $\Omega^\ast_\mu$ and an ordering on $\Omega_\mu$, that we will also denote~$\leq$. 

Let $\mathfrak{A}:=(A;\mathbf{R})$ be a relational structure with signature $\mu=(\mu_i;i\in I)$.
The {\em skeleton}, $\skel(\mathfrak{A})$, of $\mathfrak{A}$ is
the set $\{\mathfrak{A}_{\restriction F}: F \text{\,  is a finite
subset of $A$}\}$.  The {\em age}, $\age(\mathfrak{A})$, of
$\mathfrak{A}$ is the set of isomorphism  types of the elements of
$\skel(\mathfrak{A})$.  If $\mathfrak{B}$ is a relational
structure then, by a slight abuse of notation, we allow ourselves
to write $\mathfrak{B}\in \age(\mathfrak{A})$ to indicate that the
isomorphism type of $\mathfrak{B}$ is an element of the age of
$\mathfrak{A}$. Note that $\age(\mathfrak{A})$ with the relation
$\leq$ is also a poset.

A subset $\mathcal{A}$  of $\Omega_\mu$,  is an {\it ideal} if :
\begin{enumerate}
\item $\mathcal A$ is non-empty.
\item $\mathcal A$ is an {\em initial segment}, that is $\mathfrak B\in \Omega_\mu$, $\mathfrak C\in \mathcal A$ and $\mathfrak B\leq \mathfrak C$ implies  $\mathfrak B\in \mathcal A$.
\item  $\mathcal A$ is {\em up-directed}, that is $\mathfrak B, \mathfrak B'\in \mathcal A$ implies $\mathfrak B, \mathfrak B'\leq \mathfrak C$ for some $\mathfrak C\in \mathcal A$.
\end{enumerate}


Clearly, the  age  of a relational structure  is an ideal. As
shown by  Fra\"{\i}ss\'e, see \cite{fraissetr},  the converse
holds for countable ideals and,  hence,   particularly  in the case
when  $\mu$ is finite. If $I$ is infinite, the converse also holds for every ideal consisting  of the finite models  of a set of universal first-order sentences (see Section \ref{section:representableideals}).  In his book, W. Hodges proposed  to find an ideal for which  the converse does not hold as an exercise for which he has
no solution, see \cite{hodges} Exercise 17 Chapter 7,  p.332.  Such an ideal, obtained with the first author,  is  described in \cite{pou-sob}. It is made of binary relational structures  coding metric spaces which isometrically embed into the real line equipped with the ordinary distance. Because of the existence of   such  an example, we may say that an ideal $\mathcal A$ of $\Omega_{\mu}$  is {\it
representable}  if there is some relational structure $\mathfrak
A$ such that $\age(\mathfrak {A})= \mathcal A$ and it is
$\kappa$-representable   if there is some relational structure
$\mathfrak A$ of cardinality  $\kappa$  such that $\age(\mathfrak
{A})= \mathcal A$. 

One purpose of this paper is to point out  the  following
\begin{problem} \label{mainproblem}
Which ideals of $\Omega_{\mu}$ are representable, which ideals are not?
\end{problem}

We present only partial results. We give first some examples of representable ideals, see Subsection \ref{universal}.  Examples lead us to  consider the same problem for ideals made of finite metric spaces,  ordered by isometrical embedding, see Subsection  \ref{submetric}.  The special case of ideals included into $age(\R^n)$, where $\R^n$ is equipped with the euclidian distance, is left unresolved.  
In Subsection \ref {examples}, 
we provide many more examples of ideals made of binary relational structures which are not representable. They are based on a notion of {\em ashes}. 

In Theorem \ref{thm:one} below,  we will
characterize those signatures $\mu$ for which every ideal of
$\Omega_{\mu}$ is representable. 

\begin{theorem} \label{thm:one}
The following statements are equivalent: 
\begin{enumerate} [{(i)}]
\item  Every ideal of $\Omega_{\mu}$ is representable.
\item The set $\mu^{-1}[\N^*\setminus \{1\}]$ is finite.
\item The set  $ \mathscr{J}( \Omega_{\mu})$ of ideals of $\Omega_{\mu}$,  equipped with the product topology  on $\mathfrak P(\Omega_{\mu})$,  is compact.
\end{enumerate}
\end{theorem}

\begin{problem}\label{compact}
Let $\mathcal A$ be an ideal of $\Omega_\mu $.  If the set $\mathscr{J}(\mathcal A)$  of ideals included into $\mathcal A$, equipped with the product topology  on $\mathfrak P(\mathcal A)$,  is compact, does $\mathcal A$  have a representation? 
\end{problem}

Note that if  $\mathscr{J}(\mathcal A)$ is compact,  this is the Stone space of the Boolean algebra generated by the subsets of $\mathcal A$ of the form $\uparrow\hskip -2pt \mathfrak s:= \{\mathfrak t\in \mathcal A: \mathfrak s\leq \mathfrak  t\}$ for $\mathfrak s\in \mathcal A$(cf \cite{bekkali-pouzet-zhani}). 
We may represent members of this Boolean algebra by "sentences", replacing $\uparrow\hskip -2pt \mathfrak s$ by $\exists \mathfrak s$ and $\mathcal A\setminus\hskip-3pt \uparrow\hskip -2pt  \mathfrak s$ by $\forall \neg \mathfrak s$ (this can be made more concrete by means of  infinitary sentences).  A positive solution of Problem \ref{compact} above  amounts to a compactness theorem (for a counterpart, see Subsection \ref{universal}).

The other purpose of this paper is to derive from this study a solution  of a long standing question of  Cusin and Pabion \cite {cusin-pabion}. 
 
 Indeed, on  $\Omega^{*}_{\mu}$ the same  notion of ideals can be
introduced. Since $\Omega^{*}_{\mu}$ is a proper class, we extend
the above stipulations by requiring that an ideal should be a set
(and not a proper class). We say that an 
ideal $\mathcal{A}$ of $\Omega^\ast_\mu$ is {\em bounded} if all its elements have cardinality less than some cardinal $\kappa$; it is   {\em $\kappa$-bounded}  if all  elements have cardinality less
than $\kappa$ and for every $\lambda<\kappa$ there is an element
of $\mathcal{A}$ of cardinality $\lambda$. It follows that every
infinite  ideal of $\Omega_\mu$ is an $\aleph_0$-bounded ideal of
$\Omega^\ast_\mu$.

In \cite {cusin-pabion} Cusin and Pabion generalized the notions
of age and ideal of isomorphism types as follows. For a relational
structure $\mathfrak A$ and an ideal $\mathcal J$ of subsets of
$A$ they associated the set $\age_{\mathcal J}(\mathfrak A)$
consisting of isomorphism types of substructures induced by
$\mathfrak A$ on elements of $ \mathcal J$; more formally, $\age_{\mathcal J}(\mathfrak
A):=\{Is(\mathfrak{A}_{\restriction S}) \mid S\in \mathcal{J}\}$.
If isomorphism types are quasi-ordered by embeddability, this set is an
ideal of the quasi-ordered set $\Omega^\ast_\mu$.

Let us say that an ideal  $\mathcal{C}$ of $\Omega^\ast_\mu$   is
{\it representable}, if there is a relational structure $\mathfrak
A$ and an ideal $\mathcal J$ of  subsets of $A$ such that
$\age_{\mathcal J}(\mathfrak A)=\mathcal C$. Note that if
$\mathcal{C}$ is an ideal of $\Omega_\mu$ then it is representable
in this more general sense if and only if it is representable as
defined previously. (To check this, let $\mathfrak{A}$ be a
relational structure and let $\mathcal{J}$ be an ideal of finite
subsets of $A$  so that
$\age_\mathcal{J}(\mathfrak{A})=\mathcal{C}$. Let $B$ be the union
of the elements of $\mathcal{J}$ and
$\mathfrak{B}:=\mathfrak{A}_{\restriction B}$. Then every finite
subset $F$ of $B$ is in $\mathcal{J}$  because the singletons of
$B$ are elements of $\mathcal{J}$ and  $\mathcal{J}$ is an
updirected initial segment. Hence $\mathfrak{B}_{\restriction
F}\in \mathcal{C}$ implying $\age(\mathfrak{B})\subseteq
\mathcal{C}$. On the other hand every element of $\mathcal{C}$ is
finite and hence an element of $\age(\mathfrak{B})$.)

Cusin and Pabion asked the following question. Suppose that $\mu$
is a singleton signature. Is it then true, that every ideal of
$\Omega^\ast_\mu$ is representable? The answer is negative. In
fact we will prove in Theorem \ref{thm:cusin}, that if $\mu$ is a
singleton ternary signature then there is an ideal of
$\Omega^\ast_\mu$  whose elements are countable relational
structures  and $\mathcal A$ is not representable.

\begin{theorem}\label{thm:cusin}
Let $\mu$ be a singleton ternary  signature. There is an   $\aleph_1$-bounded ideal $\mathcal A$ in $\Omega^\ast_\mu$   which is not representable.
\end{theorem}

We do not know if there is an example  with binary relations. In
the case when $\mu$ is a singleton signature, we do not know if
for every uncountable cardinal $\kappa$ there exists a
non-representable $\kappa$-bounded ideal of $\Omega^\ast_\mu$.

\section{Representable and non-representable  ideals} \label{section:representableideals}

\subsection{Ideals defined by sets of universal sentences}\label{universal}
A sufficient   condition for
representability of an ideal $\mathcal A$ of  $\Omega_{\mu}$ with an  arbitrary  signature $\mu$ can be expressed in  model-theoretic terms.  As it is well-known, the class $mod(T)$ of models of a first-order theory $T$  is an ideal of $\Omega^{*}_{\mu}$ if and only if $T$ is universal (that is can be  axiomatized by universal sentences) and for every disjunction $\varphi \vee \psi$ of universal sentences, $\varphi \vee \psi \in T$ if and only if $\varphi\in T$ or $\psi\in T$ \cite{chk}.  With the compactness theorem of first -order logic, it follows  that an ideal $\mathcal A$ of  $\Omega_{\mu}$,  which consists of the finite models (up-to isomorphism) of a universal theory,  is representable. Furthermore, if $\mathcal A$ is infinite, it is $\kappa$-representable for every $\kappa\geq \mathcal \vert A\vert$. 

The  above condition on $\mathcal A$ can be easily translated in terms of \emph{reducts} as follows. 
  
  \indent Let 
$\mu: I\to \N^{*}$ be a signature. If  $I'$ is a subset of $I$, we denote by $\mu_{\restriction I'}$ the restriction of $\mu$ to $I'$. If  $\mathfrak {A}:=(A; (\mathsf{R}^{\mathfrak{A}}_i)_{i\in I})$ is a relational structure with signature  $\mu$,   the \emph{$I'$-reduct of}  $\mathfrak {A}$ is the relational structure ${\mathfrak A}^{\restriction I'}:= (A; (\mathsf{R}^{\mathfrak{A}}_i)_{i\in I'})$ of signature $\mu_{\restriction I'}$. If $I'$ is finite, ${\mathfrak A}^{\restriction I'}$ is a {\it finite reduct} of $\mathfrak {A}$.  If $\mathcal C$ is a class of relational structures of signature $\mu$, we denote by $\mathcal C^{\restriction I'}$ the class of $I'$-reducts of members of $\mathcal C$. We denote by $\widehat {\mathcal C}$ the class of relational structures $\mathfrak A$ such that $\mathfrak A^{\restriction I'}\in \mathcal C^{\restriction I'}$ for every finite $I'\subseteq I$. We say that $\mathcal C$ is \emph{closed} if $\mathcal C= \widehat { \mathcal C}$. We use freely the same notations  for classes made of isomorphism types of relational structures. 

It is not hard to show that if $\mathcal A$ is an ideal of $\Omega_{\mu}$ then $\widehat{\mathcal A}$ is an ideal too. And also that an ideal $\mathcal A$ is closed if and only if $\mathcal A$ is the set of finite models of a universal theory. Hence, we may recast the aforementioned fact as:
\begin{theorem} Every closed ideal $\mathcal A$ of $\Omega_{\mu}$ is representable; if $\mathcal A$ is infinite,  it  is $\kappa$-representable for every $\kappa\geq \vert\mathcal A\vert $.
\end{theorem}
A proof using compactness is a straightforward exercice.  See  \cite{pou-sob} for a   more detailed discussion.
\subsection{The extension property}
Let $\mathcal A$ be an ideal of $\Omega_\mu$.  A relational structure $\mathfrak A$ with age included into $\mathcal A$ is \emph{ extendable w.r.t. $\mathcal A$}  if for every $\mathfrak B\in \mathcal A$  there is some $\mathfrak C$,  with age included into $\mathcal A$,  which extends both $\mathfrak A$ and $\mathfrak B$. 
An  ideal  $\mathcal{A}$  of $\Omega_\mu$ has the {\em extension
property} if   every $\mathfrak A\in \Omega^{*}_\mu$ such that  $\vert A \vert <\kappa := \vert  \mathcal{A}\vert $ and  $age (\mathfrak A)\subseteq \mathcal{A}$ is extendable  w.r.t. $\mathcal A$. 

\begin{lemma}\label{extension}   If an ideal $\mathcal A$ of $\Omega_\mu$ has the extension property then it is  representable.
\end{lemma}
\begin{proof} Let  $\kappa:= \vert \mathcal A\vert$ and let $(\mathfrak  B_{\alpha})_{\alpha <\kappa}$ be an enumeration of the members of $\mathcal A$. We define a sequence $(\mathfrak A_{\alpha})_{\alpha<\kappa}$ such that:
\begin{enumerate}
\item $\vert A_\alpha\vert<\aleph_{0}$ if $\alpha<\omega$ and $\vert A_\alpha\vert\leq\vert  \alpha \vert $, otherwise.
\item $A_\alpha \subseteq A_{\alpha'}$ and $(\mathfrak{A}_{\alpha'})_{ \restriction A_\alpha}=\mathfrak A_{\alpha}$  for every $\alpha <\alpha^\prime <\kappa$.
\item $\mathrm{age}(\mathfrak A_{\alpha })\subseteq \mathcal A$.
\item $\mathfrak B_{\alpha}\leq \mathfrak A_{\alpha +1}$.
\end{enumerate}
We start with $\mathfrak A_{0}$ equal to the relational structure
on the empty set, and we use transfinite recursion. To get
$\mathfrak A_{\alpha+1}$ we apply the extendibility property of $\mathcal A$ to 
$\mathfrak A:= \mathfrak A_{\alpha}$ and $\mathfrak B:=
\mathfrak B_{\alpha}$. At limit stages we define $\mathfrak
A_{\alpha}$ to be $\cup_{\gamma < \alpha} {\mathfrak A}_{\gamma}$.
Clearly, $\mathfrak A_{\kappa}$ has age $\mathcal A$.
\end{proof}

\begin{corollary} \label{repideal}Every countable ideal is representable.
\end{corollary}
In view of Problem \ref{compact} we may ask:
\begin{problem} Let $\mathcal A$ be an ideal   of $\Omega_\mu $. If   $\mathscr{J}(\mathcal A)$   is compact, does $\mathcal A$ have the extension property?
\end{problem}
\subsection{The amalgamation property}
Let  $\mathcal C\subseteq \Omega^*_{\mu}$.  Let  $f_1: \mathfrak A\to \mathfrak A_1$ and $f_2: \mathfrak A\to \mathfrak A_2$ be a pair of  embeddings such that   $\mathfrak A,  \mathfrak A_1,  \mathfrak A_2\in \mathcal  C$. We say that this pair  \emph{amalgamates}  if there are two embeddings $g_1: \mathfrak A_1\to \mathfrak B$ and $g_2: \mathfrak A_2\to \mathfrak B$ such that  $\mathfrak B\in \mathcal C$ and $g_1\circ f_1=g_2\circ f_2$. We say that $\mathcal C$ has the \emph{amalgamation property} if every pair of embeddings amalgamates. If this property holds  for pairs of embeddings whose  domain have size at most $\kappa$, we say that $\mathcal C$ has the \emph{$\kappa$-amalgamation property}
\begin{lemma} Let  $\mathcal A$ be an ideal of $\Omega_{\mu}$;
if  $\mathcal A$ has the amalgamation property, then the collection of  countable $\mathfrak A$  whose age is included into $\mathcal A$ has the $\aleph_0$-amalgamation property. In particular, if $\mathcal A$ has size at most $\aleph_1$ then it has the extension property. 
\end{lemma}
\begin{proof} One proves first that every pair  of embedding  $f_1: \mathfrak A\to \mathfrak A_1$ and $f_2: \mathfrak A\to \mathfrak A_2$ such that   $\mathfrak A, \mathfrak A_1 \in \mathcal A$, $age(\mathfrak {A}_2)  \subseteq  \mathcal  A$ and $A_2$ countable amalgamates. For that, one  writes $A_2$ as an increasing sequence $(A_{2,n})_{n\in \omega}$ of finite sets containing the image of $A$ and one successively amalgamates $\mathfrak A_1$ with the $\mathfrak A_{2, n}$. This allows to do the same  when the condition  on $\mathcal A_1$ is relaxed.
\end{proof}
\begin{corollary} If an ideal of $ \Omega_{\mu}$ has the amalgamation property  and has size at most $\aleph_ 1$ then   it is representable. \end{corollary}

\begin{problems} \label{secondproblem}Let $\mathcal A$ be an ideal of $\Omega_{\mu}$. Suppose that $\mathcal A$ has the amalgamation property. 
\begin{enumerate}
\item Does $\mathcal A$ has a representation?
\item Is there an homogeneous $\mathfrak A$ such that $age(\mathfrak A)= \mathcal A$?
\end{enumerate}
\end{problems}

\subsection{Metric spaces as relational structures and representability}\label{submetric}

Metric spaces can be encoded, in several ways,  as binary relational structures in such a way that isometries correspond to embeddings.  For example, let $\mu: \Q_{+}\rightarrow \{2\}$. To each metric space $\M:= (M, d)$,  where $d$ is a distance over the set $M$,  we may associate the relational structure $rel(\M):=(M, (\delta_r)_{r\in \Q_{+}})$ of signature $\mu$ where $\delta_r(x,y)=1$ if $d(x,y)\leq r$ and $\delta_r(x,y)=0$ otherwise.  With this definition,  $d(x,y)$ is the infimum of the set of $r$'s such that $\delta_{r}(x,y)=1$, hence we may recover $\M$ from $rel(\M)$. From this fact, it follows that:
\begin{enumerate}
\item for two metric spaces $\M:= (M, d)$, $\M':=(M', d')$, a map $f: M\rightarrow M'$  is an isometry from $\M$ into $\M'$ if and only this  is an embedding from $rel(\M)$ into $rel(\M')$.  
\item Moreover, if $\mathfrak A$ is a relational structure of signature $\mu$ such that every induced substructure on at most $3$ elements embeds into $\M'$, then there is a  distance $d$ on $A$ such that $rel((A,d))= \mathfrak A$.  
\end{enumerate}
If we compare metric spaces via isometric embeddings, the class   $\mathcal M$, resp. $\mathcal M_{<\omega}$,  of  metric spaces, resp. finite metric spaces, is an ideal. Hence, $\mathcal M$ and  $\mathcal M_{<\omega}$ yield an ideal of $\Omega^{*}_{\mu}$ and of $\Omega_{\mu}$ respectively. It make sense then to consider the representability of  an ideal  $\mathcal C$ of $\mathcal M_{<\omega}$. Because of item 1 above, its  image $\mathcal A$  into  $\Omega^{*}_{\mu}$ is an ideal and because of item 2  the representability of $\mathcal A$ amounts to the representability of $\mathcal C$. 

The ideal $\mathcal M_{<\omega}$ is representable, eg by the space $\ell^{\infty}_{\R}(\N)$ of bounded sequences of reals, equipped with the "sup" distance. But, it turns out that there are plenty of non-representable ideals of $\mathcal M_{<\omega}$. We give some examples below. 

Let $\M:= (M, d)$ be  a metric space. Let $a\in M$, we set  $spec(\M, a):=\{d(a, x):  x \in \M\}$  and, for $r\in \R_{+}$, we set $B_{\M}(a, r):= \{x\in M: d(a,y)\leq r\}$. The \emph{spectrum} of $\M$ is the set $spec(\M):= \bigcup \{spec(\M, a): a\in M\}$. The \emph{diameter of $\M$} is $\delta(\M):= sup(spec(\M))$ and we set  $d(\M):= inf ( spec(\M)\setminus \{0\})$ (hence $\delta(\M):= +\infty$ if $\M$ is unbounded and $d(\M):=+\infty$ if $\vert M\vert \leq 1$).  If $\mathcal C$ is a set of metric spaces, we set $d(\mathcal C):=inf \{d(\M): \M \in \mathcal C\}$. Let $t\in \R_{+}$, we set $\omega_t(\M):=sup \{\vert X \vert  : X\subseteq  M \; \text{and} \; d(\M_{\restriction X})\geq t\}$. We say that $\M$ is  \emph{$t$-totally bounded} if $\omega_t(\M)$ is finite and  that $\M$ is \emph{totally bounded} if $\M$ is  $t$-totally-bounded for every $t\in   \R_{+}^*$. We say that $
\M$ is \emph{$t$-uniformly  bounded}  if $\omega_t(\M_{\restriction X})\leq \varphi_t(\delta (\M_{\restriction X}))$ for  some non-decreasing map $\varphi_t: \R_{+}\rightarrow \R_{+}$ and  every bounded subspace $\M_{\restriction X}$ of $\M$. 
\begin{lemma}
Let $\M$ be a $t$-totally bounded metric space.  Let $\mathcal C\subseteq age(\M)$ be an ideal such that $d(\mathcal C)\geq t$. Then  $\mathcal C$ is representable iff $\mathcal C$ is countable. 
\end{lemma}
\begin{proof}
 Suppose that $\mathcal C$ is representable. Let $\M':= (M', d')$ be a representation. 

{\bf Claim.}
$$\vert X'\vert =\omega_t(\M'_{\restriction X'}) \leq \varphi_t(\delta (\M'_{\restriction X'}))$$ for every bounded subset $X'$ of $\M'$.
 
{\bf Proof of the Claim. } Since $age(\M')\subseteq \mathcal C$, we have $d(\M'_{\restriction X'})\geq t$. The equality follows. If the inequality above does not hold then $X'$  contains a finite subset $X''$ with more than $ s:= \varphi_t(\delta (\M'_{\restriction X'}))$ elements. But then for some finite subset $X$ of $\M$ such that 
  $\M_{\restriction X}$ is isometric to $\M'_{\restriction X''}$ we have $\vert X\vert =\omega_t(\M_{\restriction X}) \leq \varphi_t(\delta (\M_{\restriction X}))\leq s$. A contradiction. 
 
 From our claim, each ball in $\M'$ is finite, hence  $\M'$ is countable. Thus $\mathcal C=age(\M')$ is countable. 
Conversely, if  $\mathcal C$ is countable then it is representable from Corollary \ref{repideal}. 

\end{proof}
\begin{proposition}\label{thmnonrepage} Let $\M$ be  an unbounded metric space whose group of isometries, $Aut (\M)$,  acts transitively on the elements of $M$. Suppose that for some $t\in \R^{*}_{+}$, $spec(\M)\cap [t, +\infty)$ is uncountable and every bounded subset of $\M$ is $t$-totally bounded, then $age_{t}(\M):= \{Is(\M_{\restriction X}):  d(\M_{\restriction X})\geq t \; \text{and} \; X \; \text{ is finite}\}$   is a  non-representable ideal of $age(\mathcal \M)$.  
\end{proposition}

\begin{proof}

{\bf Claim 1.} Let $t\in \R_{+}^{*}$, then $\M$ is $t$-uniformly bounded.

{\bf Proof of  Claim 1.} This follows from the fact that $Aut(\M)$ is transitive and  every bounded subspace is $t$-totally bounded. To see it,  fix $a\in M$. Let   $\varphi: \R_{+}\rightarrow \R_{+}$ defined by setting $\varphi_t(r):= \omega_t(\M_{\restriction B_{\M}(a, r)})$. Let $\M_{\restriction X}$ be a bounded subspace of $\M$ and let $r:= \delta(\M_{\restriction X})$. Since $Aut(\M)$ is transitive, $\M_{\restriction X}$ is isometric to $\M_{\restriction X'}$ for some subset $X'$  of $B_{\M}(a, r)$. Hence, $\omega_t(\M_{\restriction X})=\omega_t(\M_{\restriction X'})\leq \omega_t(\M_{\restriction X'})=\varphi_t(r)$. 

{\bf Claim 2.} $age_t(\M)$ is an uncountable ideal. 

{\bf Proof of Claim 2. } Since $\M$ is unbounded and $Aut(M)$ is transitive, $age_t(\M)$ is an ideal. Let $a\in M$. Since $Aut(\M)$ is transitive, $spec(\M, a)=spec(\M)$, hence $spec(\M, a)\cap [t, +\infty[$ is uncountable. Thus $age_t(\M)$ contains uncountably many $2$-element metric spaces.
\end{proof}

We consider more generally ideals made of metric spaces which omit a given set of distances. Precisely, let $\M$ be a metric space, $\kappa:= \vert M\vert $ and  let $A \subseteq \R_{+}^{*}$; we set   $age_{- A}(\M): =\{ Is(\M_{\restriction X}):  spec(\M_{\restriction X})\cap A= \emptyset \; \text{and} \; X \; \text{ is finite}\}$.  Given a type  $\mathfrak b\in age(\M)$, let $Orb(\mathfrak b, \M):= \{X\subseteq M: Is(\M_{\restriction X})=\mathfrak b\}$; given $r\in R_{+}$ and $a\in M$, we set $S_{\M}(a, r):=\{x\in M: d(a, x)=r\}$.

\begin{proposition} \label{omittingdistance}Let $\M$ be a metric space and let $\kappa:=spec(\M)$. Suppose that   $\kappa $ is infinite and that there is a cardinal $\lambda <\kappa:=spec(\M)$ such that  for every $\mathfrak s\in age(\M)$, $Orb(\mathfrak b, \M)$ contains a subset $\mathcal X_{\mathfrak b}$ of size at least $\kappa$ such that $\vert \{F\in \mathcal X_{\mathfrak b}: F\cap S_{\M}(a,r)\not = \emptyset \}\vert  \leq \lambda$ for every $a\in M, r\in \R_{+}$. Let $A \subseteq \R_{+}^{*}$ such that $\vert A\vert <\kappa$. Then  $age_{- A}(\M)$ is an ideal representable by some subspace of $\M$.
\end{proposition}
\begin{proof}
We mimick the proof of Lemma \ref{extension}. Let $(\mathfrak  b_{\alpha})_{\alpha <\kappa}$ be an enumeration of the members of $age_{- A}(\M)$. We define a  sequence $(X_{\alpha})_{\alpha<\kappa}$  of subsets of $M$ such that:
\begin{enumerate}
\item $\vert X_\alpha\vert<\aleph_{0}$ if $\alpha<\omega$ and $\vert X_\alpha\vert\leq\vert  \alpha \vert $, otherwise.
\item $X_\alpha \subseteq X_{\alpha'}$ 
\item $\mathrm{age}(\mathbb M_{\restriction X_{\alpha} })\subseteq \mathrm{age} _{- A}(\M)$.
\item $\mathfrak b_{\alpha}\in \mathrm{age}(\mathbb M_{\restriction X_{\alpha+1} })$.
\end{enumerate}
We start with $X_0:= \emptyset$. To get
$X_{\alpha+1}$ we select some   subset $F\in Orb(\mathfrak b_{\alpha}, \M)$ such that $\mathrm{age}(M_{\restriction X_{\alpha} \cup F})\in \mathrm{age} _{- A}(\M)$. If this was imposssible, then for each $F\in  \mathcal X_{\mathfrak b}\subseteq Orb(\mathfrak b_{\alpha}, \M)$ we will find $(x_F,  r_F)\in X_{\alpha}\times A$ such that $d(x_F, y_F)=r_F$ for some $y_F\in F$.  Since $\vert \mathcal X _{\mathfrak b} \vert \geq \kappa>\vert X_{\alpha}\times A\vert $, there is a subset $\mathcal X'$ of size at least $\lambda^+$  and a pair $(a, r)$ such that $(x_F, r_F)= (a,r)$ for all $F\in \mathcal X'$ but then $\vert \{F\in X_{\mathfrak b}: F\cap S_{\M}(a,r)\not = \emptyset\}\vert \geq \lambda^+$ contradicting our hypotheses on $Orb(\mathfrak b_{\alpha}, \M)$.  This allows to set $X_{\alpha +1}:= X_{\alpha}\cup X$. At limit stages we define $
X_{\alpha}$ to be $\bigcup \{X_{\gamma}: \gamma < \alpha\}$.
\end{proof}

Let $n\in \N^{*}$ and  let $\R^{n}$ be the set of $n$-tuples of reals, equipped with the euclidian  distance $d_2$. Then $(\R^{n}, d_2)$  satisfies the  hypotheses of Proposition~\ref {thmnonrepage} above and this for every $t$. It also satisfies the hypotheses of Proposition~\ref{omittingdistance} (fix a direction in $\R^{n}$ and in each $Orb(\mathfrak b, \R^{n})$ select an orbit  $\mathcal X_{\mathfrak b}$ according to this group and finally set $\lambda: =\aleph_{0}$). The same facts hold  if the euclidian distance is replaced by any distance associated with a vector space  norm  on $\R^n$. Then, we have the following:
\begin{corollary}\label{nonrepagereal}
Let  $n\in \N^*$. 
\begin{itemize}
\item For every positive real $t$, the set $age_t((\R^{n}, d_2))$ of  isometric types of  finite  subspaces $\X$  of $(\R^{n}, d_2)$  such that $d(\X)\geq t$ is a  non representable ideal.  
\item For every subset $A\subset \R^{*}_{+}$ of size $\kappa<2^{\aleph_0}$, the set $age_{- A}(\M)$ of isometric types of subspaces $\X$  of $(\R^{n}, d_2)$ whose distances does not belong to $A$ is an ideal representable by a subset of $\R^{n}$. 
\end{itemize}
\end{corollary}

The example of a non-representable ideal given in \cite {pou-sob} is $age_1(\R)$. By taking $A:= \Q_{+}$, the second item of the corollary above asserts that there are subspaces $\X$ of the real line whose age is the set $age_{- \Q_{+}}(\R)$  made of all finite metric spaces with no rational non-zero distances. Such spaces are sections of the quotient $\R/\Q$ of the additive group $\R$ by the additive group $\Q$, but not every section provides such a space. The metric space $\mathbb{S}$ made of the unit circle with the arc  length metric satisfies the hypotheses of Proposition~\ref{omittingdistance}. \emph{We do not know if $\mathrm{age} (\mathbb{S})$ contains an non-representable ideal}.

In the case of $\R$ or even $\R^n$, we can say a little more. 
Let  $(\R, d)$ where  $d(x,y):= \vert x-y\vert$. Note first that a $3$-element metric space isometrically embeds into $(\R, d)$ iff one distance is the sum of the two others; a $4$-element metric space whose all $3$-element subsets embed into $(\R,d)$ does not necessarily embed into $(\R,d)$ (think of four vertices forming a  "rectangle" whose sides have length $a$ and $b$ and diagonal length $c:= a+b$). However, all the $\leq 4$-element subspaces of a metric space $\M$ embed  isometrically into $(\R, d)$ if and only if $\M$ isometrically embeds into $(\R, d)$;  moreover an embedding from $\M$ into $\R$ is determined by its values on any $2$-element subset of $\M$. This extends:  \emph{all $\leq n+3$-elements subspaces of a metric space $\M$ embed into $(\R^{n}, d_2)$ iff   $\M$ embeds into $(\R^n, d_2)$}. 
From this, we immediately have:
  \begin{itemize} 
\item  if $\mathcal C \subseteq  age((\R^n, d_2))$ is  a representable  ideal, all its representations embed into $(\R^{n}, d_2)$. Hence have cardinality at most the continuum. 
  \end{itemize}
This is a substantial difference with  to the representability of closed ideals. 
 
Let us mention that in the case of the real line, the two problems in Problems 4 have a positive answer.
\begin{lemma} Let $\R$ be the real line equipped with the ordinary distance. Let $\mathcal C \subseteq  age(\R)$ be an ideal. 
Then $\mathcal C$ has the $2$-amalgamation property if and only if  there is a  homogeneous metric space $\D$ whose age is $\mathcal C$. Moreover, if  $\mathcal C$ contains at least a $3$-element metric space, then $\D=(G, d_{\restriction G})$ where $G$ is an additive subgroup of $\R$ and $d_{\restriction G}$ is induced by the distance on $R$.   
\end{lemma}

\begin{proof} 
We just give a hint.  Let  $\mathcal C \subseteq  age(\R)$ be an ideal. Let $V:=\cup \{spec(\M): \M\in \mathcal C\}$. 

{\bf Case 1.}
 $\vert V\vert \leq 2$. In this case  $V= \{0, v\}$ and  $\D:= (V, d_{\restriction V})$ has the required property.  
 
 {\bf Case 2. } $\vert V\vert \geq 2$.  Set $G:= V\cup -V$.
  
 {\bf Claim 1.} If $\mathcal C$ has the $2$-amalgamation property then $G$ is a subgroup of $\R$.
 The proof of this claim breaks into three parts; we leave the verification to the  reader. 

 Subclaim 1.  For every finite subset $F$ of $V$ there is some $\M\in \mathcal C$ and $a\in M$  such that 
 $F\subseteq spec(\M,a)$. 

Subclaim 2. $V$ is unbounded. 

Subclaim 3.   $y-x\in V$  
and $x+y\in V$ for every   $x, y\in V$ with $x\leq y$. 

{\bf Claim 2.} Let $G$ be an additive subgroup of $\R$ and $\D:= (G,d_{\restriction G})$. Let $f$ be an  isometry from a subset $A$ of $G$ onto a subset $A$. Then  $f$ extends to an isometry. 

Indeed, we may suppose $A\not =\emptyset$. Let $x\in A$ and $x':= f(x)\in A'$. Let  $g^{+}(y):= x'+y-x$ for all $y\in G$ and $g^{-}(y):= x'-y+x$ for all $y\in G$. These two maps are isometries from $\D$ into iself and one of these  extends $f$. This proves Claim~2.
\end{proof}

Note that the  $1$-amalgamation property provides a representative $\M$ with $Aut(\M)$ transitive. As an example let $M:= a\cdot\Z \cup (b+ a\cdot\Z)$ with $0<b<\frac{a}{2}$ and $\M:= (M, d_{\restriction M})$. 
 
 \begin{problems}\label{thirdproblem}
 \begin{enumerate}
 \item Describe the amagamable ideals and the homogeneous subspaces of $(\R^{n}, d_2)$;
 \item Characterize the representable ideals of $(\R^{n}, d_2)$. 
 \end{enumerate}
 \end{problems}

 
\subsection{A construction  of non-representable ideals}\label{examples}

Let $S$ be a set. A set $\mathcal{S}$   of finite subsets of  $S$  is an {\em ash on $S$} if:
\begin{enumerate}
\item $\{s\}\in \mathcal{S}$ for every $s\in S$.
\item For every finite subset $\mathcal{\mathcal F}$ of $\mathcal{S}$ there exists an element $s\in S\setminus\bigcup\mathcal{F}$ so that $\{s\}\cup F\in \mathcal{S}$ for every set $F\in \mathcal{\mathcal F}$.
\item For every subset $S'$ of $S$ with $|S'|=|S|$  there is a finite subset $F$ of $S'$ with $F\not\in \mathcal{S}$.
\end{enumerate}
Note that there is no ash on a finite set.

An ash can be obtained as follows. Let $\kappa$ be an infinite cardinal and $\mathcal{T}$ a family which consists of $\kappa$ sets of    cardinality $\kappa^+$ for which there is no finite subset $F\subseteq \bigcup\mathcal {T}$ with  $X\cap F\not=\emptyset$ for all $X\in \mathcal {T}$. (For example the elements of $\mathcal {T}$ are disjoint.) Let $0\not=n\in \omega$ and let $\mathcal{S}$ be the set of finite subsets $F\subseteq S:=\bigcup{\mathcal {T}}$ with $|F\cap X|\leq n$ for all $X\in \mathcal{T}$. Then $\mathcal{S}$ is an ash on $S$.

Let  $\mu:= S\rightarrow \{2\}$.  Let $\mathcal{S}$ be an ash on  $S$ and let $\mathcal {A_{\mathcal S}}$ be the collection of all  finite relational structures $\mathfrak {A}$ in $\Omega_{\mu}$ for which for all $x,y,z\in A$ and $s,t\in S$:
\begin{enumerate}
\item $\neg R_s(x,x)$.
\item If $R_s(x,y)$ then $R_s(y,x)$.
\item If $R_s(x,y)$ and $R_t(x,y)$ then $s=t$ and if $R_s(x,y)$ and $R_s(x,z)$ then $y=z$.
\item If  $x\not=y$ then there exists an element $r\in S$ so that $R_r(x,y)$.
\item Every non-empty subset of the set $\{r\in S : \exists u\in A : R_r(x,u) \}$ is an element of $\mathcal{S}$.
\end{enumerate}
Note that the  elements of $\mathcal {A_{\mathcal S}}$ are graphs with several types of edges, a type of edge for every element of $S$.

\begin{lemma}\label{lem:ages}
Let $\mathcal{S}$ be an ash on the set $S$. Then $\mathcal {A_{\mathcal S}}$   is an ideal of $\Omega_{\mu}$.
\end{lemma}
\begin{proof}
It follows directly from the definition that $\mathcal {A_{\mathcal S}}$ is closed under induced substructures. Let $\mathfrak{A}$ and $\mathfrak{B}$ be two elements of $\mathcal {A_{\mathcal S}}$ with $A\cap B=\emptyset$. Item 2 of the definition of ash allows us to determine successively an edge type for every pair $(x,y)$ with $x\in A$ and $y\in B$, satisfying items 1 to  5 of the definition of $\mathcal {A_{\mathcal S}}$. It follows that $\mathcal {A_{\mathcal S}}$ is updirected.\end{proof}

\begin{lemma}\label{lem:notrep}
Let $\mathcal{S}$ be an ash on the set $S$. Then $\mathcal {A_{\mathcal S}}$ is not representable. 
\end{lemma}
\begin{proof}
Assume for a contradiction that there is a relational structure  $\mathfrak {A}$ whose age is equal to $\mathcal {A_{\mathcal S}}$.

Let $x\in A$. It follows from  items 3 and 4 of the definition of $\mathcal {A_{\mathcal S}}$   that there exists an injection  $f: A\setminus\{x\}\to S$ so that $R_{f(y)}(x,y)$ for every element $y\in A\setminus\{x\}$.

Because every  two element structure of $\mathcal {A_{\mathcal S}}$ is isomorphic to an induced substructure of $\mathfrak{A}$ it follows that $|A|\geq |S|$.    Hence $|f[A\setminus\{x\}]|=|S|$ which in turn implies using item 3 of the definition of an ash   that there is a finite subset $F\subseteq f[A\setminus\{x\}]$  with $F\not\in \mathcal{S}$. But this leads to a contradiction because the substructure of $\mathfrak{A}$ induced by the set $\{x\}\cup \{y : \exists s\in F:  R_s(x,y)\}$ is not an element of $\mathcal {A_{\mathcal S}}$ according to item 5 of the definition of $\mathcal {A_{\mathcal S}}$.\end{proof}

With the theorem of R. Fra\"\i ss\'e asserting that every ideal  with countable signature is representable, this yields: 
\begin{corollary}\label{cor:frai}
Let $\mathcal{S}$ be an ash on the set $S$. Then $|S|>\aleph_0$.
\end{corollary}

The graph  $\mathrm{G}$ is an {\em ash-graph} if it has  the following two properties:
\begin{itemize}
\item For every finite subset $F\subseteq G$ there exists a vertex $v\in G\setminus F$ which is adjacent to every vertex in $F$.
\item The graph $\mathrm{G}$ does not contain a complete subgraph $\mathrm{K}$ with $|K|=|G|$.
\end{itemize}
Note that the set of finite subsets $F$ of $G$ which contain an element adjacent to all the other elements of $F$ is an ash on $G$.

An ash-graph $\mathrm{G}$ can be obtained as follows. Let $\M:= (M;d)$ be a metric space with the properties:
\begin{itemize}
\item For every finite subset $F$ of $M$ there exists an element $x\in M$ with $d(x,y)\geq 1$ for all $y\in F$.
\item Every subset $W$ of $M$ with $|W|=|M|$ contains two elements $x,y$ with $x\not= y$ and $d(x,y)<1$.
\end{itemize}
Such a metric space is an {\em ash-space}. The set of real numbers is an example of an ash-space.

Let $\M:= (M;d)$ be an ash-space. Then the graph  $\mathrm{G}$  with vertex set $M$ in wich two different vertices are adjacent if and only if their distance is larger than or equal to 1 is an ash-graph. On the other hand if $\mathrm{G}$ is an ash-graph, then the metric space $(G;d)$ with $d(x,y)=1$ if $x$ is adjacent to $y$ and $d(x,y)={1\over 2}$ if $x$ is not adjacent to $y$  is an ash-space.

Let $\mathrm{P}:=(P;\leq)$ be an up-directed poset  which does not contain a maximal element and no chain of size $|P|$. Such a poset is an {\em ash-poset}. Let $\mathrm{P}:=(P;\leq)$ be an ash-poset. Let $\mathcal{P}$ be the set of finite subsets $F$ of $P$ which contain an element $x$ with $x\geq y$ for all $y\in F$. It follows that $\mathcal{P}$ is an ash on $P$.

For example, let $\kappa>\aleph_0$ be a cardinal and let $\mathrm{P}$ be the poset  on the set of finite subsets of $\kappa$ with $\subseteq$ as the order relation. Then $\mathrm{P}$ is an ash-poset. 

Let $\mathcal{S}$ be an ash on the set $S$.  We used a particular construction to obtain a non-representable age. There are many other ways. For an example,  we can generalise in an obvious way from binary to $n$-ary relations and we can define generalized edges of some type $s\in S$ as $n$-tuples for which a relation of the form $R_s$ holds.  Of course the relations do not have to be symmetric.
%



%

\section{Proof of Theorem \ref{thm:one}}

We will use the following fact. Let $P$ be a poset and let $\mathscr{J}(P)$ be the set of ideals of $P$. We think of $\mathscr{J}(P)$ as being equipped with the topology induced by the product topology on the power set of $P$.   Then \emph{ $\mathscr{J}(P)$  is compact if and only if $P$ is a finite union of principal final segments and 
$\uparrow\hskip -2pt \hskip-4pt x\, \, \cap \uparrow\hskip -2pt \hskip -4pt y$ is a   finite union of principal final segments,  for all $x,y\in P$ } (for a proof, see \cite{bekkali-pouzet-zhani}).

Applying this to a non empty  initial segment $\mathcal C$ of $\Omega_\mu$ and observing that $\Omega_\mu$ has a least element and does not have an infinite descending chain,  we get that $\mathscr{J}(\mathcal C)$ is compact if and only if  for every $\mathfrak A, \mathfrak B\in \mathcal C$  there are at most  finitely many  non-isomorphic $\mathfrak C\in \mathcal C$ such that:
\begin{enumerate}
\item $\mathfrak A, \mathfrak B\leq \mathfrak C$.
\item If  $x\in \mathfrak C$ then  $\mathfrak A\not \leq \mathfrak C_{\restriction C\setminus \{x\}}$ or  $\mathfrak B\not \leq \mathfrak C_{\restriction C\setminus \{x\}}$.
\end{enumerate}

\begin{lemma}\label{lemcompact} $\mathscr{J}(\Omega_\mu)$ is compact if and only if $\mu^{-1}[\mathbb{N^*}\setminus\{1\}]$ is finite.
\end{lemma}
\begin{proof}
Let us check that if $\mathcal C:= \Omega_{\mu}$ and $\mu^{-1}[\mathbb{N}\setminus\{1\}]$ is finite then only finitely many non-isomorphic $\mathfrak C$ satisfy  Conditions 1 and 2 above.
Let $\mathfrak A, \mathfrak B\in \mathcal C$. Let $m:= \vert
A\vert$ and $n:= \vert B\vert $. Suppose $\mathfrak C$ satisfies
 Conditions $1$ and $2$  and let $r_{\mathfrak C}:= \vert C\vert $. We may suppose
that $C=\{1,\dots, r_{\mathfrak C}\}$. Let $A_{\mathfrak C},
B_{\mathfrak C}\subseteq C$ such that $\mathfrak C_{\restriction
A_{\mathfrak C}}\cong \mathfrak A$ and  $\mathfrak C_{\restriction
B_{\mathfrak C}}\cong \mathfrak B$.  Clearly $A_{\mathfrak C}\cup
B_{\mathfrak C}= C$, hence $r_{\mathfrak C}\leq m+n$. If there are
infinitely many non-isomorphic $\mathfrak C$ satisfying Conditions $1$ and $2$,
there are infinitely many for which $r_{\mathfrak C}$,
$A_{\mathfrak C}$, $B_{\mathfrak C}$, $\mathfrak C_{\restriction
A_{\mathfrak C}}$ and  $\mathfrak C_{\restriction B_{\mathfrak
C}}$ are independent of $\mathfrak C$. Let $r,  A,  B$ such a triple. Since
$A\cup B= \{1, \dots, r\}$, all unary relations on  $\{1, \dots,
r\}$ are entirely determined. Let $I':= \mu^{-1}[\mathbb{N^*}\setminus\{1\}]$. Since $I'$  is finite, the number of relational structures of signature  $\mu_{\restriction I'}$ defined on $\{1, \dots, r\}$ is finite but then one cannot define infinitely many relational
structures of signature  $\mu$ on this set. A contradiction.

Let $\underline 2$ be the constant map from $\N$ to $\{2\}$. 

{\bf  Claim 1. } If  $\mu^{-1}[\mathbb{N^*}\setminus\{1\}]$ is  infinite then $\mathscr{J}(\Omega_{\underline 2})$ can be mapped continuously into $\mathscr{J}(\Omega_{\mu})$ by a one-to-one map.

{\bf Proof of Claim 1. } Let $\varphi: \omega \rightarrow \mu^{-1}[\mathbb{N^*}\setminus\{1\}]$ be a one-to-one map and let  $rang(\varphi)$ be its range. For every  $\mathfrak A:= (A, (R_{i})_{i<\omega})\in \Omega^*_{\underline 2}$, let $F(\mathfrak A):= (A, (S_i)_{i\in I})$ where $S_i: A^{n_i}\rightarrow \{0\}$ if $i\not \in rang(\varphi)$ and $S_i((x_{1},\dots, x_{n_i})):=R_j ((x_1,x_2))$ if  $i:=\varphi (j)$. Clearly,
\begin{enumerate}
\item $F(\mathfrak A_{\restriction B})= F(\mathfrak A)_{\restriction B}$ for every   $\mathfrak A\in \Omega^{*}_{\underline 2}$ and $B\subseteq A$;

\item If $\mathfrak A$, $\mathfrak A'\in \Omega^{*}_{\underline 2}$,  a map $f$ is an isomorphism from $F(\mathfrak A)$ onto $F(\mathfrak A')$if and only if  $f$ is an isomorphism  from  $\mathfrak A$ onto $\mathfrak A'$.
\end{enumerate}

Consequently, $F$ defines an embedding from $\Omega_{\underline 2}$ onto an initial segment of $\Omega_\mu$. This map induces a continuous embedding from $\mathscr{J}(\Omega_{\underline 2})$  into $\mathscr{J}(\Omega_{\mu})$.

{\bf Claim 2.} $\mathscr{J}(\Omega_{\underline 2})$ is not compact.

{\bf Proof of Claim 2. } Let $\mathfrak A, \mathfrak B$, where $A:= \{0\}$, $\mathsf{R}^\mathfrak{A}_i:A^2\rightarrow \{1\}$ ,  $B:= \{1\}$, $\mathsf{R}^\mathfrak{B}_i:B^2\rightarrow \{0\}$ for $i<\N$. Let $\mathfrak C_n$ where $C_n:=\{0,1\}$, $\mathsf{R}^\mathfrak{C}_n(x,y)=1$ iff $(x,y)\in \{(0,0), (0,1)\}$ and $\mathsf{R}^\mathfrak{C}_i(x,y)=1$ iff $(x,y)=(0,0)$ in case $i\not = n$. The $\mathfrak C_n$ 's satisfy Conditions 1 and 2 above. Hence, $\mathscr{J}(\Omega_{\underline 2})$ cannot be compact.

It follows from Claims 1 and 2, that $\mathscr{J}(\Omega_{\mu})$ cannot be
compact if  $\mu^{-1}[\mathbb{N^*}\setminus\{1\}]$ is infinite. This completes the proof.
\end{proof}

\begin{lemma}\label{truc} $\mu^{-1}[\mathbb{N^*}\setminus\{1\}]$ is finite if and only if  every ideal $\mathcal A \in \mathscr{J}(\Omega_{\mu})$ is representable.
\end{lemma}
\begin{proof} Suppose that $\mu^{-1}[\mathbb{N^*}\setminus\{1\}]$ is finite. Let  $\mathcal A\in \mathscr{J}(\Omega_{\mu})$.  For each $\mathfrak s\in \mathcal A$, let $\uparrow\hskip -2pt  \mathfrak s:=\{\mathfrak t \in \mathcal A: \mathfrak s\leq \mathfrak t\}$. The set $\mathcal F:= \{X\subseteq \mathcal A:  \uparrow\hskip -2pt  \mathfrak s\subseteq X \mbox {\;  for \; some\; } \mathfrak s\in \mathcal A\}$ is a  filter. Let $\mathcal U$ be an ultrafilter on $\mathcal A$ containing it.
Let $\mathcal A^*$ be the set of  non-empty members of  $\mathcal A^*$. For each $\mathfrak s\in \mathcal A^*$,  let $\mathfrak S_{\mathfrak s}$ such that $Is(\mathfrak S_{\mathfrak s})=\mathfrak s$ and  $S_{\mathfrak s}:= \{1, \dots, \vert
\mathfrak s\vert \}$. Let $\mathfrak A:= \Pi_{\mathfrak s\in \mathcal A^{*} } \mathfrak
S_{\mathfrak s} / {\mathcal U}$ be the ultraproduct of the $\mathfrak S_{\mathfrak s}$' s .
Let $B:= \{(x_{\mathfrak s})_{\mathfrak s\in\mathcal A^{*}}:  \mbox {there\; is\; a }\mathfrak  t\in
\mathcal A \; \text {such that}\;  \{\mathfrak s:   \mathfrak S_{\mathfrak s \restriction \{x_\mathfrak {s}\}} =\mathfrak t \}\in \mathcal U\}$. 

{\bf Claim. } $\age(\mathfrak A_{\restriction B})=\mathcal A$.

 First, $\mathcal A\subseteq \age (\mathfrak
A_{\restriction B})$. Indeed, let $\mathfrak s\in \mathcal A$. For every $\mathfrak s'\geq \mathfrak s$
select an embedding $\varphi_{\mathfrak s'}$  of $\mathfrak S_{\mathfrak s}$ into $ \mathfrak S_{\mathfrak s'}$. Let $X:=
\{(x^{i}_{\mathfrak s'})_{\mathfrak s'\in \mathcal A}:  1\leq i\leq   \vert \mathfrak s\vert \; \text{and}  \; x^{i}_{\mathfrak s'}:= \varphi_{\mathfrak s' }(i)$ for each $\mathfrak s'\geq \mathfrak s\}$. Then
$X\subseteq B$ and $\age( \mathfrak A_{\restriction X})=\mathfrak s$. 

Next, $\age (\mathfrak A_{\restriction
B})\subseteq \mathcal A$.  Indeed, let $Y:= \{(y^{i}_{\mathfrak s})_{\mathfrak s\in \mathcal
A}: 1\leq i\leq n \}$ be a $n$-element subset of $B$. 
We claim that here is some $\mathfrak s\in \mathcal A$ such that the
projection $p_{\mathfrak s}$ from $B$ onto  $S_{\mathfrak s}$ induces an isomorphism from
$\mathfrak A_{\restriction Y}$ onto $\mathfrak S_{\mathfrak s  \restriction
\{y^{i}_{\mathfrak s}: 1\leq i\leq  n \}}$. Due  to the choice of the ultrafilter, it is  obvious  that there is some $\mathfrak s$ such that $p_{\mathfrak s}$ preseves the unary relations. Now, let $I': =\mu^{-1}[\mathbb{N^*}\setminus\{1\}]$; since there are only finitely many relational structures of signature $\mu_{\restriction I'}$ on an $n$-element set, we may find $\mathfrak s$ such that the other relations can be preserved. From this $\mathfrak A_{\restriction Y}\in \mathcal A$. 

Suppose that $\mu^{-1}[\mathbb{N^*}\setminus\{1\}]$ is infinite.  According to Claim 1 of Lemma \ref{lemcompact},  $\mathscr{J}(\Omega_{\underline 2})$ can be mapped continuously into $\mathscr{J}(\Omega_{\mu})$ by a one-to-one map. According to Corollary \ref{nonrepagereal}, $\Omega_{\underline 2}$ contains non-representable ideals. If $\mathcal A$ is a non-representable ideal of $\Omega_{\underline 2}$ then, as it is easy to check,  its image is a non-representable ideal of $\Omega_{\mu}$. 

With this,  the proof of Lemma \ref{truc} is complete. \end{proof}

\section {Proof of Theorem \ref{thm:cusin}}
Let $\underline 2:  \N \rightarrow \{2\}$ and let $\Omega^{*}_{(3)}$ be the class of relational structures containing a single ternary relation.  To   each  relational structure  $\mathfrak{A}:=(A; \mathbf{R})\in \Omega^{*}_{\underline 2}$ such that $A\cap \N= \emptyset $,  we associate $F(\mathfrak{A}):=
(A\cup \N, T)\in \Omega^{*}_{(3)}$ such that  $T:=X\cup Y$,  where $X:=  \{(x,y,z)\in \N^3\setminus \{(1,0,1)\}: x+y=z\}\cup \{(1,0,0)\}$ and $Y:=  \{(x,y,z)\in A^2\times \N: (x,y)\in \mathsf{R}^{\mathfrak{A}}_z\}$.

{\bf Claim 1.} Let  $\mathfrak{A}:=(A; \mathbf{R})$, $\mathfrak{A'}:=(A'; \mathbf{R'})$ be two binary relational structures as above, then:
\begin{enumerate}
\item A map $f: A\rightarrow A' $ is an isomorphism from $\mathfrak{A}$ into $\mathfrak{A'}$ if and only if the map $F(f):= f\cup 1_{\N}$ extending $f$ by the identity on $\N$ is an isomorphism from  $F(\mathfrak{A})$ into $F(\mathfrak{A'})$.
\item Every isomorphism $g:F(\mathfrak{A})\rightarrow F(\mathfrak{A'})$ is of the form $F(f)$ for some isomorphism $f: \mathfrak{A}\rightarrow \mathfrak{A'}$.
\end{enumerate}
{\bf Proof of  Claim 1. } Part $(1)$ is straightforward to check.

Part $(2)$. Let $g:F(\mathfrak{A})\rightarrow F(\mathfrak{A'})$ be an isomorphism. First,  $g$ is the identity on $\N$. Indeed,  as  is is easy to see,  each element of $\N$ is definable by an existential  formula. To be precise,   $0$ is the unique element $x$ of $A'\cup \N$ such that $(x,x,x)\in T'$, hence $g(0)=0$. Also, $(0,y,z)\in T'$ implies that $y=z$ and $y\in \N$, from which follows that $g(\N)\subseteq \N$. Furthermore, $g(1)= 1$ since $1$ is the only element $x\not = 0$ of $A'$ such that $(x,0,0)\in T'$. Since $n+1$ is the unique element $x$ of $A'$ such that
$(n,1, x)\in T'$, we have $g(n+1)=n+1$ for all $n$. Second, since $g$ maps $\N$ onto $\N$, it maps $A$ into $A'$. Clearly, for every $z\in \N$, we have $(x,y)\in  \mathsf{R}^{\mathfrak{A}}_z$ if and only if $(g(x), g(y))\in  \mathsf{R}^{\mathfrak{A'}}_z$ that is $f:= g_{\restriction A}$ is an isomorphism from $\mathfrak{A}$ into $\mathfrak{A'}$, proving that $g=F(f)$. \endproof

Let $\mathcal C$ be a set of relational structures, we denote by $\downarrow\hskip -2pt  \mathcal C$ the set of isomorphism types of  relational structures which embed into some member of $\mathcal C$.

{\bf Claim 2. } Let $\mathcal C$ be a  subset of $\Omega^{*}_{\underline 2}$ made of  relational structures as above. Then $\downarrow\hskip -2pt  \mathcal C$ is a representable ideal if and only if $\downarrow\hskip -2pt  F[\mathcal C]$ is a representable ideal of $\Omega^{*}_{(3)}$.

{\bf Proof of Claim 2. } 
Suppose that $\downarrow\hskip -2pt  \mathcal C$ is representable. Let $\mathfrak A:= (A; \mathbf{R})$ be a relational structure and  $\mathfrak J$ be an ideal of   subsets of $A$ such that  $\mathcal A_{\mathfrak J}(\mathfrak A)=\mathcal \downarrow\hskip -2pt  C$. With no loss of generality,  we may suppose $A\cap \N= \emptyset$.  Let $\mathfrak B:= F(\mathfrak A)$ and $\mathfrak K:= \{X\subseteq A\cup \N: X\cap A \in \mathfrak J\}$. From Part  $(1)$ of Claim $1$  we obtain that   $\mathcal A_{\mathfrak K}(\mathfrak B)=\downarrow\hskip -2pt  F[\mathcal C]$ proving that $\downarrow\hskip -2pt  F[\mathcal C]$ is representable.

Conversely, suppose that $\downarrow\hskip -2pt  F[\mathcal C]$ is representable. Let $\mathfrak B:= (B; T)$ be a relational structure consisting of a single ternary relation $T$ and  let $\mathfrak K$ be an ideal of   subsets of $B$ such that  $\mathcal A_{\mathfrak K}(\mathfrak B)=\downarrow\hskip -2pt  F[\mathcal C]$.  For each $\mathfrak A\in \mathcal C$, there is some  $X_{\mathfrak {A}}\in \mathfrak K$ and an  isomorphism  $g_{\mathfrak{A}}$  from $F(\mathfrak{A})$ onto $\mathfrak B_{\restriction X_{\mathfrak {A}}}$.  We claim that $g_{\mathfrak {A} \restriction \N}$ is  independent of $\mathfrak {A}$.  To see  that,  take  $\mathfrak A,  \mathfrak A'\in \mathcal C$.  Since $\mathfrak K$ is an ideal, it contains the union 
$Y$ of the ranges of  $g_{\mathfrak{A}}$ and $g_{\mathfrak{A'}}$ hence there is some $Y'\in \mathfrak K$ such that  $\mathfrak B_{\restriction Y'}$ is isomorphic to a member of $F[\mathcal C]$ and there is some embedding $h$ from $\mathfrak B_{\restriction Y}$  into $\mathfrak B_{\restriction Y'}$. According to Part $(2)$ of Claim $(1)$, $h\circ  g_{\mathfrak{A}}= h\circ  g_{\mathfrak{A'}}$ proving our claim. Identifying  $\N':= g_{\mathfrak{A}}[\N]$ to $\N$ allows us to define a relational structure $\mathfrak A$ such that  $\mathfrak B= F(\mathfrak A)$.  For $\mathfrak J:= \{X\setminus \N:  X\in \mathfrak K\}$ we have
$\mathcal A_{\mathfrak J}(\mathfrak A)=\downarrow\hskip -2pt  \mathcal C$  proving that $\downarrow\hskip -2pt  \mathcal C$ is representable.

Taking for $\mathcal C$ a non-representable ideal of finite binary relational structures (as given by Corollary \ref {nonrepagereal}), we get an ideal of  countable ternary relations which is not-representable.
With this the proof of Theorem \ref{thm:cusin} is complete. \endproof

\section{Conclusion}
We just scratched the surface of Problem \ref{mainproblem}.  We posed  the question of the representability of ideals of metric spaces. Besides Problems  \ref{compact}, \ref{secondproblem}, and \ref{thirdproblem}, we offer  a very basic  one: 
\begin{problem} \label{lastproblem}
 Does the representability of an ideal of $\Omega_\mu$ depend only upon its order structure? 
In the most general form, the problem is this. Let $\mathcal A, \mathcal A'$ be two ideals of $\Omega_\mu, \Omega_{\mu'}$ which are order isomorphic; is $\mathcal A$ representable if and only if $\mathcal A'$ is representable?
\end{problem}

\end{document}